\documentclass[10pt, twoside, a4paper]{article}
\usepackage{graphicx, graphics,latexsym,amssymb, amsmath}
\usepackage[mathscr]{euscript}
\usepackage{url}
\usepackage[english]{babel}

\title{
 On triangle meshes\\
  with valence $6$ dominant vertices}
\author{Jean-Marie Morvan}
\date{}
\newcommand{\Addresses}{{
  \bigskip
  \footnotesize

  Jean-Marie Morvan,
\textsc{Universit\'{e} de Lyon, CNRS UMR $5208$,
Universit\'{e} Claude Bernard Lyon $1$, Institut Camille Jordan,
 $43$ blvd du $11$ Novembre $1918$, $F$-$69622$ Villeurbanne-Cedex, France}
 {\rm and}
 \textsc{King Abdullah University of Science and Technology, $C.E.M.S.E$,  Thuwal $23955$-$6900$, Saudi Arabia}\\
\\\par\nopagebreak
  \textit{E-mail address}~: \texttt{morvan@math.univ-lyon1.fr} {\rm and}
\texttt{Jean-Marie.Morvan@KAUST.EDU.SA}
}}

\def\cal{\mathcal}
\font\bb=msbm10
\def\Z{\hbox{\bb Z}}
\def\R{\hbox{\bb R}}
\def\E{\hbox{\bb E}}

\def\N{\hbox{\bb N}}

\def\C{\hbox{\bb C}}

\newtheorem{lemma}{Lemma}
\newtheorem{theorem}{Theorem}
\newtheorem{proposition}{Proposition}
\newtheorem{definition}{Definition}
\newtheorem{corollary}{Corollary}

\begin{document}

\maketitle

\begin{abstract}{\footnotesize}
We study triangulations $\cal T$ defined on a closed disc $X$ satisfying the following condition~: In the interior of $X$, the valence of all vertices of $\cal T$ except one of them (the irregular vertex) is $6$. By using a flat singular Riemannian metric adapted to $\cal T$, we prove a uniqueness theorem when the valence of the irregular vertex is not a multiple of $6$.  Moreover, we exhibit  non isomorphic triangulations on $X$ with the same boundary, and with a unique irregular vertex whose valence is $6k$, for some $k >1$.
\end{abstract}

\section{Introduction}
The goal of this paper is to study topological properties of triangulations on surfaces with Riemannian tools. More precisely, we study triangulations of a topological disc satisfying the following property~: Every interior vertex except one has valence $6$. This situation has been explored in \cite{li2010editing}, and in  \cite{peng2011connectivity},  \cite{peng2013connectivity} (with quad meshes instead of triangulations). Our results can be considered as their generalisation.

\section{Basic definitions and notations}
\subsection{Weighted marked points on a surface}
 Let $X$ be a topological surface $X$.  We denote by $\stackrel{\circ}{X}$ its interior and by $\partial X$ its possible boundary.
A {\it weighted marked point} in $X$ is a couple $(v, {\rm val}_v)$, where $v$ is a point of $X$ and ${\rm val}_v$ is an integer (called the {\it valence} of $v$).  If
$${\bf val}_{\stackrel{\circ}{X}}=((v_1, {\rm val}_{v_1}) ..., (v_p, {\rm val}_{v_p}))$$ is a (finite) sequence of  weighted points in $\stackrel{\circ}{X}$, and
$${\bf val}_{\partial X}=((w_1, {\rm val}_{w_1}) ..., (w_q, {\rm val}_{w_q}))$$ is a (finite) sequence of (ordered) weighted points in $\partial X$,  we denote by
$$(X, {\bf val}_{\stackrel{\circ}{X}}, {\bf val}_{\partial X})$$
 the data of $X$,  its interior with its weighted marked points and its boundary with its weighted marked points. In this case, $X$ is called a {\it weighted marked surface}. If $\stackrel{\circ}{X}$ admits  a unique weighted point $(v, {\rm val}_{v})$, we put
 $$(X, (v, {\rm val}_{v}),{\bf val}_{\partial X}) = (X, {\bf val}_{\stackrel{\circ}{X}}, {\bf val}_{\partial X}).$$

\subsection{Adapted triangulations}
We begin with a classical definition of a triangulation, useful in our context and related to the theory of {\it Dessins d'Enfants}.
\begin{definition} \cite{grothendieck1984esquisse} \cite{voisin1980cartes} \cite{malgoire12cartes}
A {\rm triangulation}  is a couple $(X, {\cal T})$, where $X$ is an (oriented) topological  surface with or without boundary $\partial X$ and ${\cal T}$  a  finite graph on it such that $X \backslash {\cal T}$ is a finite union of disjoint topological discs whose boundary are triangles (that is, the union of three edges of ${\cal T}$).
\end{definition}

\noindent Let $(X, {\bf val}_{\stackrel{\circ}{X}}, {\bf val}_{\partial X})$ be a a weighted marked surface, with $${\bf val}_{\stackrel{\circ}{X}}=((v_1, {\rm val}_{v_1}) ..., (v_p, {\rm val}_{v_p})) \mbox{ and }
{\bf val}_{\partial X}=((w_1, {\rm val}_{w_1}) ..., (w_q, {\rm val}_{w_q})).$$

\begin{definition}
 A triangulation $\cal T$ defined on  $X$ is  {\rm adapted} to $({\bf val}_{\stackrel{\circ}{X}}, {\bf val}_{\partial X})$ if
 \begin{itemize}
 \item
 its vertices lying on the boundary of $X$ are $\{w_1,...,w_q\}$ with respective valences $\{{\rm val}_{w_1},..., {\rm val}_{w_q}\}$,

 \item
 $\{v_1,..., v_p \}$ is a subset of the set of interior vertices of $\cal T$, with respective valences $\{{\rm val}_{v_1}, ..., {\rm val}_{v_p}\}$,

 \item
 and if other interior vertices have valence $6$.
 \end{itemize}
\end{definition}
\noindent (As usual, the {\it valence} of a vertex $v$ of a triangulation is the number of edges incident to $v$.)

\begin{definition}
 Let ${\cal T}_1$ and ${\cal T}_2$ be two triangulations of $X$, such that  ${\cal T}_1$ and ${\cal T}_2$ are adapted to $({\bf val}_{\stackrel{\circ}{X}}, {\bf val}_{\partial X})$.
\begin{itemize}
   \item
We say that ${\cal T}_1$  is {\rm equivalent} to  ${\cal T}_2$ if there exists a homeomorphism of  $X$ inducing a isomorphism (of graph) between the graphs ${\cal T_1}$ and ${\cal T}_2$.
    \item
We denote by  ${\cal T}_{{\bf val}_{\stackrel{\circ}{X}}, {\bf val}_{\partial X}}$ the set of equivalence classes of such triangulations on $X$.
\end{itemize}
\end{definition}

\subsection{Triangulations with a unique "irregular" vertex}
In the following, we will focus on triangulations whose all interior vertices have valence $6$, except one of them. This leads to the following definition.

\begin{definition}
  A triangulation $\cal T$ of a surface $X$ (with or without boundary) is said  {\rm of type} $(6,n)$ if it satisfies the following property~: The valence of each vertex interior to $X$ is $6$, except exactly one vertex (whose valence is $n \neq 6$).
\end{definition}

\noindent The goal of this paper is to prove the following result~:

\begin{theorem} \label{TTHH}
Let $(X, (v, {\rm val}_{v}),{\bf val}_{\partial X})$ be a weighted marked closed topological disc, and  $n \in \N^*$.
\begin{enumerate}
\item \label{TTHH1}
If $n$ is not a multiple of $6$,  there exists (up to an isomorphism), at most one triangulation of type $(6,n)$ adapted to $(X, (v, {\rm val}_{v}),{\bf val}_{\partial X})$.

\item \label{TTTH2}
One can build a weighted marked closed topological disc $(X, (v, {\rm val}_{v}),{\bf val}_{\partial X})$ endowed with non isomorphic triangulations of type $(6,n)$, where $n$ is a multiple of $6$.
\end{enumerate}
\end{theorem}

\noindent The proof (given in Section \ref{PROO}) of Theorem \ref{TTHH}-\ref{TTHH1} consists of enlarging the framework by introducing the class of Riemannian flat metrics with conical singularities on $X$, and use a developability argument. For Theorem \ref{TTHH}-\ref{TTTH2}, we build an explicit example, by using a covering argument.

\section{Flat metric with conical singularities on a surface}
 \subsection{On Riemannian metrics with conical singularities}
We give here classical definitions and results concerning Riemannian metrics with conical singularities on a surface, \cite{troyanov1986surfaces}, \cite{troyanov1991prescribing}, \cite{troyanov2007moduli}.

\begin{definition} \label{DEA}
Let $X$ be a  surface, endowed with marked points
$$\{v_1, \ldots, v_p\} \subset \stackrel{\circ}{X}, \{w_1, \ldots, w_q\} \subset \partial X.$$
Let $X'=X \backslash \{v_1, \ldots, v_p, w_1, \ldots, w_q\}$. A metric $g$ {\rm with conical singularities} at $\{v_1, \ldots, v_p,  w_1, \ldots, w_q\}$ on $X$ is a Riemannian metric on $X'$  such that
each  point $v_i$ ({\it resp.} $w_i$) admits  a neighborhood with polar coordinates $(r, \phi_i)$, where $r$ denotes the distance to $v_i$  ({\it resp.} $w_i$) and $\phi_i \in \R/\theta_i$ is the angular variable, $\theta_i$ being the angle of $v_i$  ({\it resp.} $w_i$).
\end{definition}

 \noindent Then, each marked point admits a neighborhood isometric to a Euclidean cone with angle $\theta_i$. A more precise and complete definition using an analytic point of view is the following. We denote by $(x,y)$ the coordinates of the Euclidean plane $\E^2$ and identify $\E^2$ with $\C$. If $z \in \C$ we write  $z=x+iy$.

\begin{definition} \label{DEB}
Let $g$ be a (singular) Riemannian metric defined on $X$.
\begin{enumerate}
\item
A point $v$ belonging to the interior of $X$ is a {\rm conical} singularity of order
$\mbox{\rm ord}_g(v)=\beta >-1$ (and of angle $\theta= 2\pi(\beta + 1)$) for $g$ if there exists a chart
  $$\psi : U \to V \subset \E^2,$$
  where $U$ is an open subset of $X$ containing $v$,
such that
$$\psi(v)=0 \mbox{ and } g \underset{U}{=}(e^{2u}|z|^{2\beta}|dz|^2);$$

\item
A point $w$ belonging to  $\partial X$ is a {\rm conical} singularity for $g$ of order
$\mbox{\rm ord}_g(w)=\gamma >-\frac{1}{2}$ (and of angle $\theta= 2\pi(\gamma + \frac{1}{2})$) if there exists a  chart
$$\psi : U \to V \subset \{(x,y) \in \E^2, y \geq 0\},$$
where $U$ is an open subset of $X$ containing $v$,
such that $$\varphi(w)=0 \mbox{ and } g \underset{U}{=} e^{2u}|z|^{4\gamma}|dz|^2.$$

\noindent In these two items,  $u$ is a continuous function, of class $C^2$ on $U \backslash \{v\}$, ({\rm resp.} $U \backslash \{w\}$). Therefore, if $v$ ({\rm resp.} $w$) is an interior ({\rm resp.} boundary) singular point,  $U$ is isometric to a neighborhood of the vertex of a cone of angle $\theta$ in the Euclidean space $\E^3$ ({\it resp.} a neighborhood of the vertex of an angular sector of angle $\theta$ in the Euclidean space $\E^2$).

 \item
 If $\{v_i, w_j\}_{i\in I, j \in J}$ is a finite set of  conical singularities of $g$ of respective orders $\{\beta_i, \gamma_j\}_{i\in I, j \in J}$,
 the formal sum ${\bf \beta} = \sum \beta_{i}v_{i} + \sum \gamma_{j}w_{j}$ is called the {\rm divisor} representing $g$. The {\rm support} of ${\bf \beta}$ is the set $\{v_i, w_j\}_{i\in I, j \in J}$. The {\rm degree} of ${\bf \beta}$ is $|\beta|=\sum \beta_i + \sum \gamma_j$.
 \end{enumerate}
\end{definition}

\noindent The link between Definition \ref{DEA} and Definition \ref{DEB} is the following~: The metric  can be written as follows~:
\begin{itemize}
\item
Around each singular interior  point,
  \begin{equation}
  g \underset{U}{=} dr^2 + r^2d\varphi^2 \underset{U}{=} e^{2u}|z|^{2\beta}|dz|^2,
  \end{equation}
\item
Around each singular boundary  point,
  \begin{equation}
  g \underset{U}{=} dr^2 + r^2d\varphi^2 \underset{U}{=} e^{2u}|z|^{4\beta}|dz|^2.
  \end{equation}
\end{itemize}

\begin{definition}
Let $X$ be a (compact) Riemannian surface with divisor $\beta$, and topological Euler characteristic $\chi(X)$. The real number
$$\chi(X, {\bf \beta})= \chi(X) + |{\bf \beta}|$$ is called the {\rm Euler characteristic} of $(X, {\bf \beta})$.
\end{definition}

\noindent Gauss-Bonnet Theorem can be stated as follows~:
\begin{theorem} \label{GAUBO}
Let $(X,g)$ be a (compact) Riemannian surface with conical singularities, curvature $K$ in its interior, and whose boundary $\partial X$ has geodesic curvature $k$. Then,
$$\int_X K da + \frac{1}{2\pi}\int_{\partial X}kds=\chi(X, {\bf \beta}).$$
\end{theorem}

\noindent As a direct consequence of Theorem \ref{GAUBO}, we get~:

\begin{corollary}\label{TODI}
Let $(X,g)$ be a  closed topological disc endowed with a Riemannian flat metric with conical singularities.  Then,
$$ \frac{1}{2\pi}\int_{\partial X}kds - \sum_{w_j \in \partial X} \gamma_j = 1 + \sum_{v_i \in \stackrel{\circ}{X}}\beta_i.$$
In particular,
\begin{itemize}
\item
if  $\partial X$ has a null geodesic curvature, then
$$- \sum_{v_i \in \stackrel{\circ}{X}}\beta_i -\sum_{w_j \in \partial X} \gamma_j =1.$$
\item
If  $\partial X$ has a null geodesic curvature and $g$ admits a unique conical singularity $v_0$ in $\stackrel{\circ}{X}$, then
$$\beta_0 =-1 - \sum_{w_i \in \partial X}\gamma_j.$$
In other words, in this special case, $\beta_0$ is completely determined by the conical singularities of the geodesic boundary.
\end{itemize}
\end{corollary}

\noindent {\small {\bf Proof of Corollary \ref{TODI}~-}} Indeed, if $X$ is flat with geodesic boundary, then $\chi(X, {\bf \beta})=0$. Moreover, if $X$ is a topological disc, then $\chi(X)=1$. Therefore, Corollary \ref{TODI} is proved.

\subsection{Flat discs with a unique conical singularity in its interior}
Suppose that $X$ is a topological (compact) disc endowed with a flat Riemannian metric without any singularity in $ \stackrel{\circ}{X}$. Then $X$ is developable into the Euclidean plane $\E^2$~: There exists an isometric immersion (unique, up to a rigid motion in $\E^2$)  of $X$ into $\E^2$  (in general, this immersion is not injective).
  We will now study the case where $X$ admits a finite set of conical singularities, such that all of them except one   belong to $\partial X$.  One has the following uniqueness result~:

\begin{proposition}\label{TTHH2}
Let $X$ is a (compact) topological disc. Let $g_1$ ({\it resp.} $g_2$) be a Riemannian  flat  metrics defined on $X$ with a possible finite set of  conical singularities. Suppose that
\begin{enumerate}
\item
 all singularities of $g_1$ ({\it resp.} $g_2$)  belong to $\partial X$,  except one of them (denoted by  $v_1$ ({\rm resp.} $v_2$));

\item
$g_1$ and $g_2$ coincide on $\partial X$~: In particular,  $g_{1_{\partial X}} = g_{2_{\partial X}}$, and $g_1$ and $g_2$ admit the same conical singularities $\{w_j\}_{j \in J} \subset \partial X$, with for all $w_j $, ${\rm \gamma}_{g_1}(w_j)={\rm \gamma}_{g_2}(w_j)$.
\end{enumerate}
Then,
\begin{itemize}
\item
${\rm ord}_{g_1} (v_1) = {\rm ord}_{g_2} (v_2)$,

\item
$v_1= v_2$.
\end{itemize}
\end{proposition}

\noindent Proposition \ref{TTHH2} can be interpreted as follows~: Under the assumptions of  Corollary \ref{TODI}  the  order $\beta_0$ of the unique interior conical singularity $v$ is completely determined by the conical singularities of the geodesic boundary. Proposition \ref{TTHH2} shows moreover that under the assumption on the angle of $v$, its position  is also determined.

\vspace{10pt}

\noindent Proposition \ref{TTHH2} is the consequence of the classical theory of developable surfaces and the following trivial lemma~:

\begin{lemma}\label{ISOC}
Let $a$ and $b$ be points in the Euclidean plane $\E^2$ and let $\theta$ be a real number.
\begin{enumerate}
\item
If $\theta \neq 2k\pi, k \in \Z$ and $a \neq b$, there exist exactly two points $m$ in the plane such that the triangle $amb$ is isosceles. These points are on the bisector of the segment $ab$.

\item
If $\theta=2k\pi, k \in \Z$, for any point $m$ in the plane,  $amb$ is a degenerate (isosceles) triangle.
\end{enumerate}
\end{lemma}

\noindent {\small {\bf Proof of Proposition \ref{TTHH2}~-}
 \begin{itemize}
 \item
 Since $g_1$ and $g_2$ have the same singularities on the boundary, with the same order, Corollary \ref{TODI} implies that ${\rm ord}g_1(v_1)={\rm ord}g_2(v_2)$.

 \item
 Let $g$ be a Riemannian  flat  metric defined on $X$. We suppose that  $g$ admits a unique conical singularity  $p$ in $\stackrel{\circ}{X}$, with angle $\theta \neq 2k\pi, k \geq 1$. Let $\gamma$ be any (smooth) simple curve whose interior belongs to $\stackrel{\circ}{X}$, joining $p$ to $\partial X$. Let us cut $X$ along $\gamma$.  We just built a new flat (for the metric induced by $g$) surface without singularity in its interior
$$\tilde X = X \cup \gamma \cup \gamma',$$
where we denote by $\gamma'$ the curve $\gamma$ with opposite orientation.
 Since $\tilde{X}$ is flat, simply connected and its interior has no singularity,  we can  develop $\tilde{X}$ into $\E^2$~: There exists an isometry
   $$i :\tilde{X} \to \E^2.$$
   In particular, $i(\partial \tilde{X})= i(\partial X \cup \gamma \cup \gamma')$ is a curve uniquely determined up to a rigid motion in $\E^2$. To simplify the notations, we put $i(v)=v$.
  The point $p$ gives rise to two points $p_1$ and $p_2$ in $\E^2$. Moreover, because the angle of $x$ is different to $2k\pi, k\geq 1$, $p_1 \neq p_2$. Considering the isoceles triangle $p_1xp_2$, Lemma \ref{ISOC} implies that the position of $x$ in $\E^2$ and then in $\tilde{X}$ and $X$ is uniquely determined up to a symmetry with respect to the line $p_1p_2$. We can choose $p$ far enough to $x$ to be sure that only one position of $x$ is possible. Then the position of $x$ in $\tilde{X}$ and $X$ is uniquely determined.
  \end{itemize}
}

\noindent {\bf Remark~-} The reader can be easily convinced that the assumptions of the type of singularities on the boundary can be much weakened, as far as a Gauss-Bonnet Theorem is satisfied $\ldots$
%

\section{Euclidean structure defined on a triangulation}
 If $(X, \cal T)$ is a topological surface endowed with a triangulation, one can define on $X$ a {\it flat Riemannian metrics} (with possible conical singularities at the vertices of ${\cal T}$) as follows~: One associates to each edge a length $l$ in such a way that the triangular inequality is satisfied on each triangle. Then, one associates to each face $t$ of ${\cal T}$ the metric of a triangle drawn in the Euclidean plane whose edges have the same lengths as the ones of $t$. Such a metric is called a {\it Euclidean structure} on $(X, \cal T)$. This construction endows $X$  with a flat metric with conical singularities.

\vspace{10pt}

\noindent A crucial example is obtained by  associating the length $1$ to each edge of ${\cal T}$. Such a triangulation is then called an {\it equilateral} triangulation. In this case, $X$ is endowed with a flat metric $g$ with conical singularities at some vertices  of $\cal T$. An interior vertex is singular if its valence $\mbox{\rm val}(v)$ is different to $6$, and a vertex $w$ belonging to the boundary of $X$ is singular if its valence $\mbox{\rm val}(w)$ is different to $3$. In such cases,
$\mbox{\rm ord}_g(v)= 2\pi - \mbox{\rm val}(v)\frac{2\pi}{6}$
and
$\mbox{\rm ord}_g(w)= \pi - \mbox{\rm val}(v)\frac{\pi}{3}$.


\section{Proof of Theorem \ref{TTHH}} \label{PROO}
We can now prove Theorem \ref{TTHH} by considering the canonical equilateral structure associated to $(X, {\cal T})$.

\subsection{The case $n \neq 6k$}
We endow $(X, {\cal T})$ with the structure of {\it equilateral} triangulation, by affecting the length $1$ to each edge. Then, $X$ is endowed with a structure of flat Riemannian disc with conical singularities. The result is then a direct consequence of Proposition \ref{TTHH2}.

\subsection{The case $n=6k$}
We suppose now that the valence of the unique singular vertex $v_0$ is a multiple of $6$~: $n=6k, k \in \N, k \neq 1$.

 \subsubsection{An example}
\begin{figure}[htbp]
\begin{center}
    \includegraphics[width=1.3\textwidth]{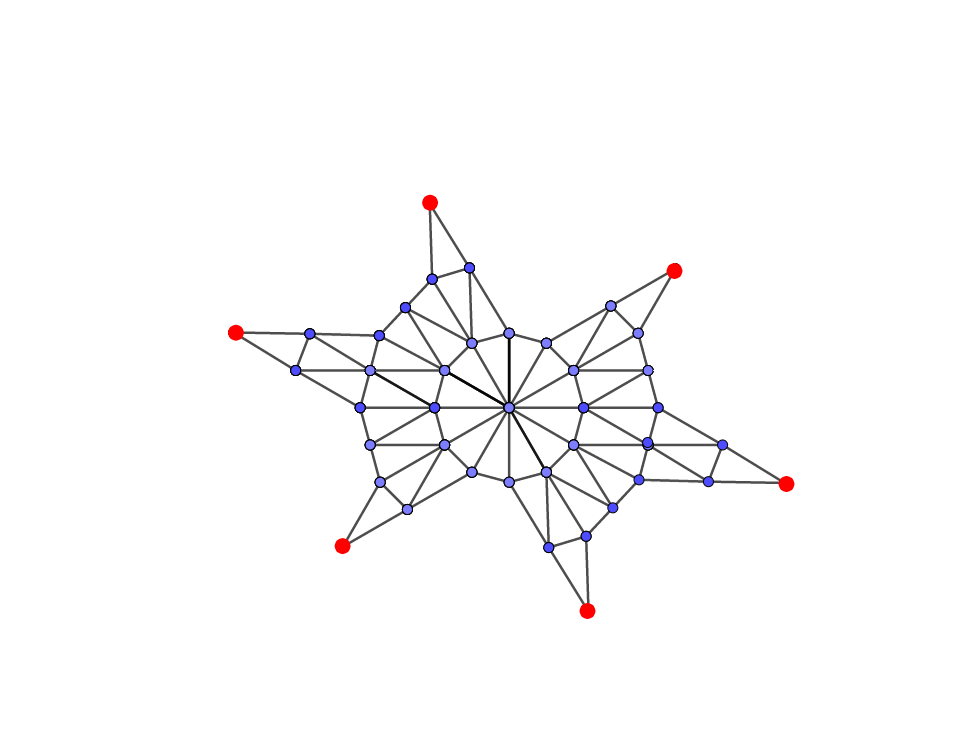}
\end{center}
     \caption{The smallest distance from the singular vertex to a corner of the boundary is $3$}
     \label{TRI25}
\end{figure}

\begin{figure}[htbp]
\begin{center}
    \includegraphics[width=1.3\textwidth]{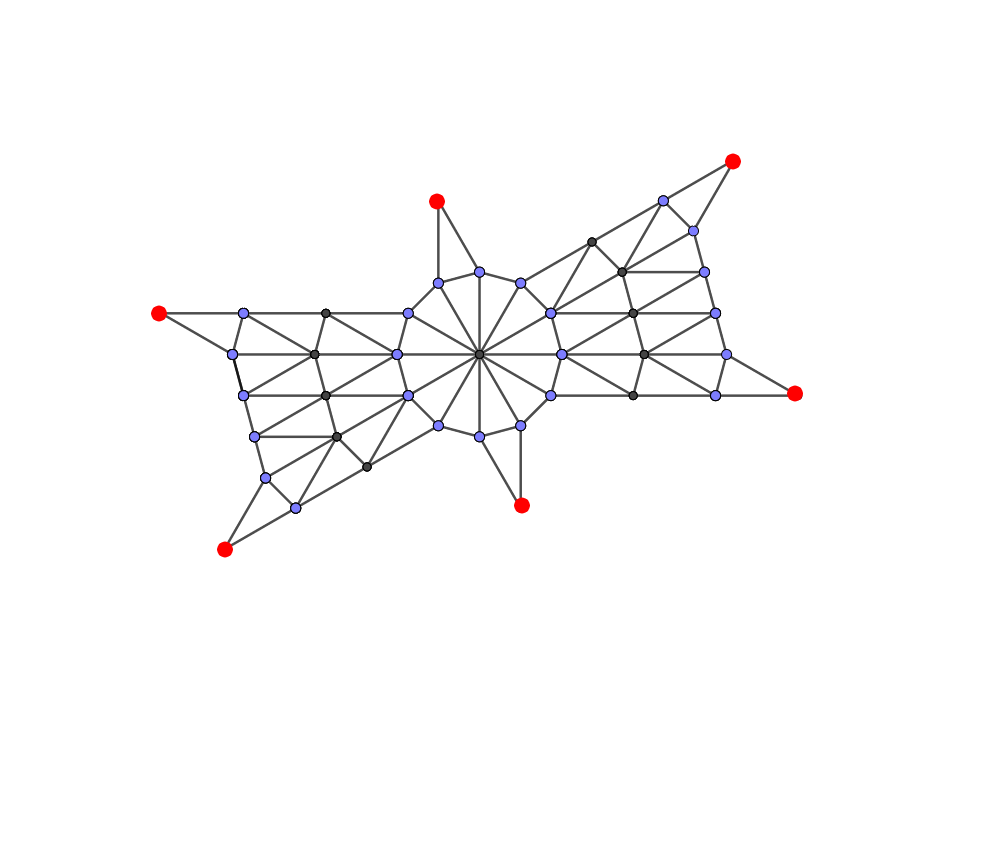}
\end{center}
     \caption{The smallest distance from the singular vertex to a corner is $2$}
     \label{TRI257}
\end{figure}

 To simplify our construction, we  assume that $k=2$ (see subsection \ref{GENE} for generalisations), and we show in Figure \ref{TRI25} and \ref{TRI257} an example of two  triangulations ${\cal T}$ and ${\cal T}'$ of type $(6,12)$ defined  on the same  (connected simply connected) domain $X$,  adapted to the same weighted marked point $(v_0, 12)$ in the interior of $X$ (the triangulations have a unique vertex of degree $12$ in $\stackrel{\circ}{X}$) and the same marked points on  ${\bf val}_{\partial X}$, that are not isomorphic.
Indeed, let us introduce the classical distance $d_{\cal T}$ on a graph~:
If $(X, {\cal T})$ is a surface endowed with a triangulation, and ${\cal V}$ is the set of vertices of $X$, the distance  $d_{\cal T}$  on ${\cal V}$ is defined as follows~: If $v_1 \in {\cal V}$ and $v_2 \in {\cal V}$,
$d_{\cal T}(v_1 , v_2)$ is the minimum number of edges connecting  $v_1$ and  $v_2$.
 In our examples,  the smallest distance $d_{{\cal T}}$ from $m$ to a corner (that is, a boundary vertex with valence $2$) of $\partial X$ is different to the smallest distance $d_{{\cal T}'}$ from $m'$ to a corner of $\partial X'$. We deduce that these triangulations are not isomorphic.

%

\subsubsection{How to build such examples}
Our construction is based on a ramified covering of a domain of the plane with two sheets, around the singular vertex.

 \subsubsection{The ramified covering of a equilateral triangulation around a conical singularity}
 Let us consider  a regular equilateral triangulation of a domain $(X, {\cal T})$ of $\E^2$ (the valence of each interior vertex is $6$), and two copies $(X_1, {\cal T}_1)$ and $(X_2, {\cal T}_2)$ of it.
Let $v_0$  be a vertex of ${\cal T}$, $v_1$ ({\it resp.} $v_2$) the corresponding vertices  in ${\cal T}_1$, ({\it resp.} ${\cal T}_2$) through the isometry between ${\cal T}$, ${\cal T}_1$ and ${\cal T}_2$. We suppose that the smallest distance of $v_0$ to a corner of ${\cal T}$ is  $2$ (with respect to  $d_{{\cal T}}$), and then, the smallest distance of $v_1$ ({\it resp.} $v_2$) to $\partial {\cal T}_1$ ({\it resp.} $\partial {\cal T}_2$) is $2$, with respect to $d_{{\cal T}_1}$ ({\it resp.} $d_{{\cal T}_2}$). Let $v_1-p_1$ be a minimal path  linking $v_1$ to  $\partial {\cal T}_1$, and $v_2-p_2$ be the corresponding minimal path in ${\cal T}_2$  linking $v_2$ to  $\partial {\cal T}_2$. We cut along these paths, so that  $v_1 -p_1$ is now replaced by two paths $v_1-p_1$ and $v_1-p'_1$ ($v_2- p_2$ is also replaced by $v_2-p_2$ and $v_2-p'_2$). Now, we identify $v_1$ and $v_2$, and we glue $v_1-p_1$ with $v_2-p'_2$, and $v_1-p'_1$ with  $v_2-p_2$. We obtain  a new domain $Y$ whose boundary is the union of $\partial {\cal T}_1$ and $\partial {\cal T}_2$ endowed with a triangulation ${\cal T}$ whose all interior vertices have valence $6$ except $v_1=v_2$ whose valence is $12$.

\vspace{10pt}

\noindent The result is a two sheets  covering $\tilde{X}$ of $X$, branched at ${\tilde v_0}$. The domain $\tilde{X}$ is endowed with an equilateral triangulation. The valence of  ${\tilde v_0}$  is $12$.

\vspace{10pt}

\noindent To build the second example, the process is exactly the same,  replacing the vertices $v_1$ and $v_2$ with minimal distance $d_{\tau_1}$ ({\it resp.} $d_{\tau_2}$) to the corners equal to $2$, by vertices $x_1$ and $x_2$ with minimal distance $3$ to the corners. We obtain  a new domain $Z$ whose boundary is the union of $\partial {\cal T}_1$ and $\partial {\cal T}_2$ endowed with a triangulation ${\cal T}'$ whose all vertices have valence $6$ except $x_1=x_2$ whose valence is $12$. Moreover, by construction, these two triangulations are not isomorphic.

\vspace{10pt}

\noindent Finally, we remark that these triangulations do not lie in the plane (they are immersed in $\E^3$ and can be embedded in $\E^4$). To draw them on the plane, one sends each triangle on the plane but not isometrically~: By remarking that a point running on  the boundary of the domain    comes back to its initial position after turning twice around $m$, we send each triangle in the plane by dividing one of its angle by two, so that the boundary of the drawing is now a simple closed curve around the image of $m$.



\subsubsection{Generalisations} \label{GENE}
\begin{figure}[htbp]
\begin{center}
    \includegraphics[width=0.7\textwidth]{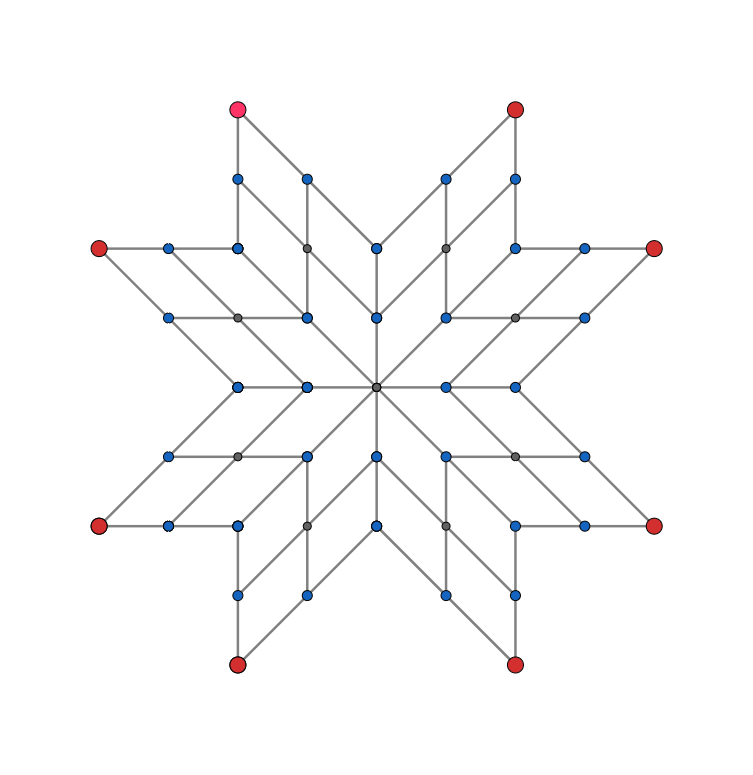}
\end{center}
     \caption{The distance of the singular vertex to the boundary is $2$ .}
     \label{QUASMESH}
\end{figure}

\begin{itemize}
\item
Analogous examples can be built by considering $k$-sheets ramified covering  giving rise to triangulations around a singular vertex of valence $6k$, for any $k >2$.

\item
Our results and constructions with triangulations can be mimicked by replacing triangles with quad meshes whose  vertices except one have valence $4$. We present here two figures showing two non isomorphic meshes of this type, see Figure \ref{QUASMESH}.

\item
The reader can be easily convinced that the assumptions on the singularities of $\partial X$ could be considerably enlightened.
\end{itemize}

\newpage

\begin{figure}[ht!]
\begin{center}
    \includegraphics[width=0.6\textwidth]{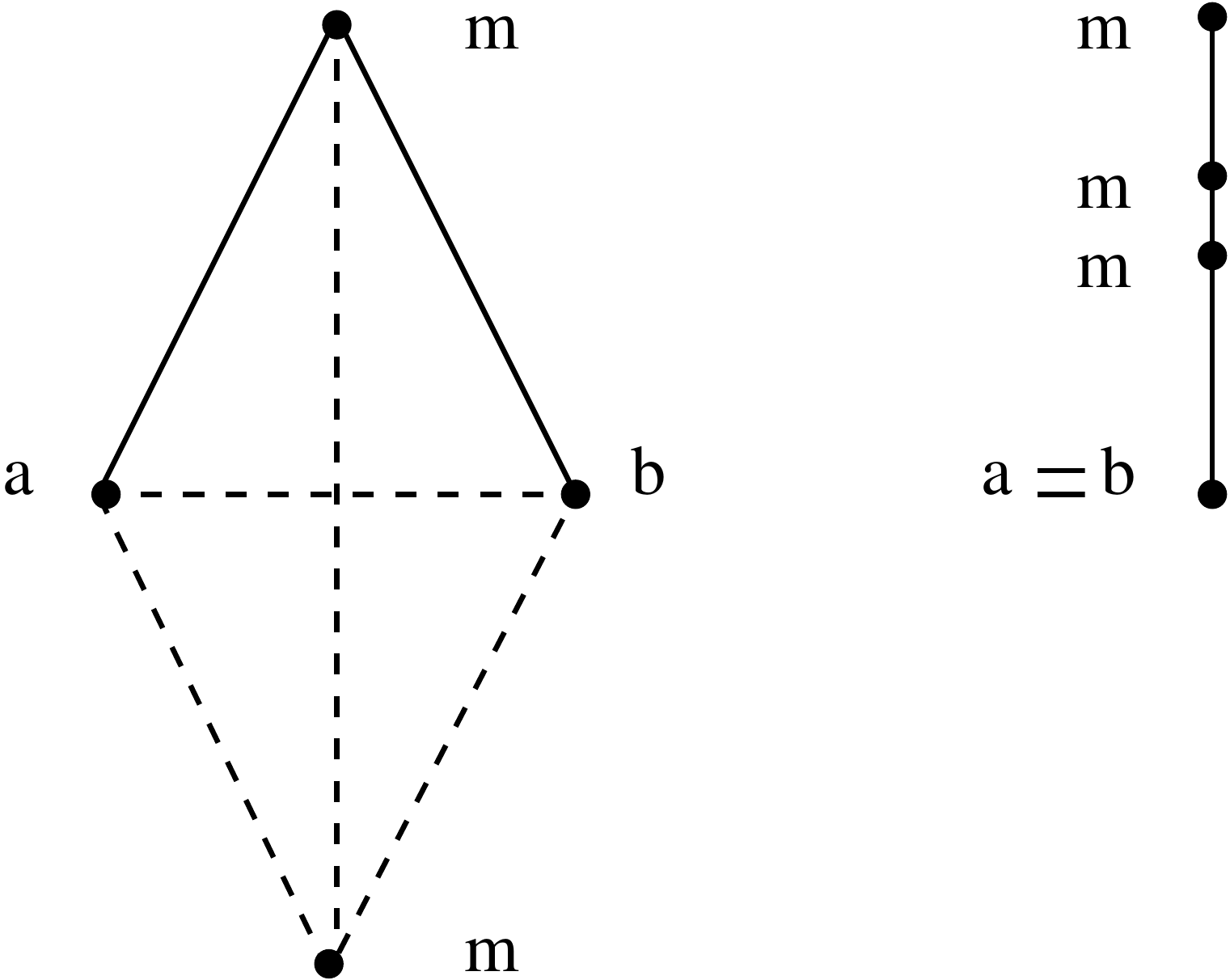}
\end{center}
     \caption{On the left, the angle $\theta$ at $m$ is not a multiple of $2\pi$. Up to a symmetry,  the point $m$ is uniquely determined by $\theta$. It is no more true on the right, where the angle at $m$ is a multiple of $2\pi$.}
     \label{ISOC2}
\end{figure}

\begin{figure}[htbp]
\begin{center}
    \includegraphics[width=0.3\textwidth]{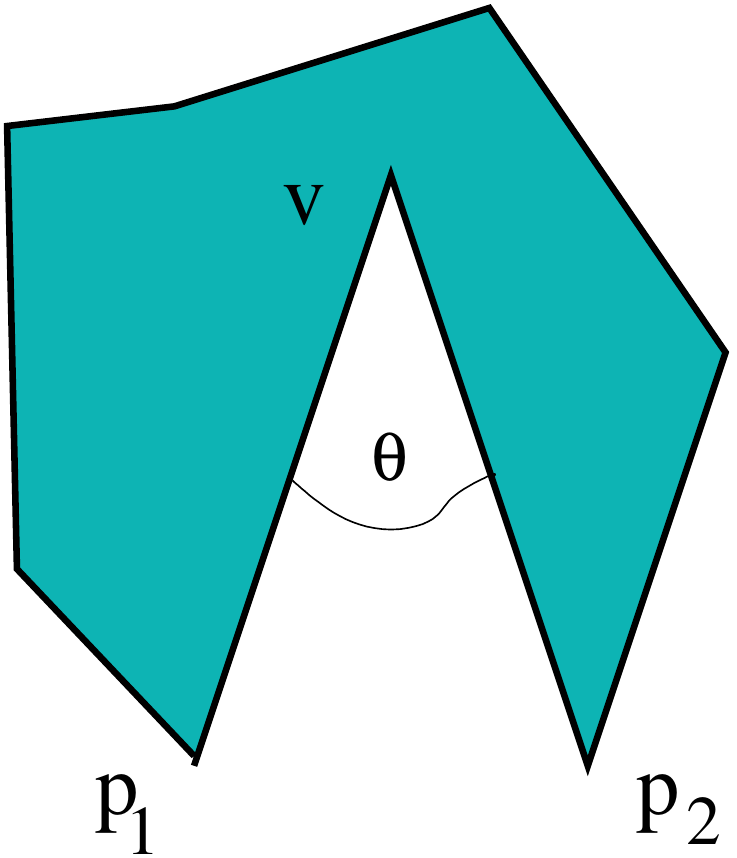}
\end{center}
     \caption{If one knows the position of $p_1$ and $p_2$ in the plane and the angle $\theta$, one knows the position of $v$ in $\tilde{X}$ and then in $X$.}
     \label{ISOC2}
\end{figure}

\begin{figure}[htbp]
\begin{center}
    \includegraphics[width=0.9\textwidth]{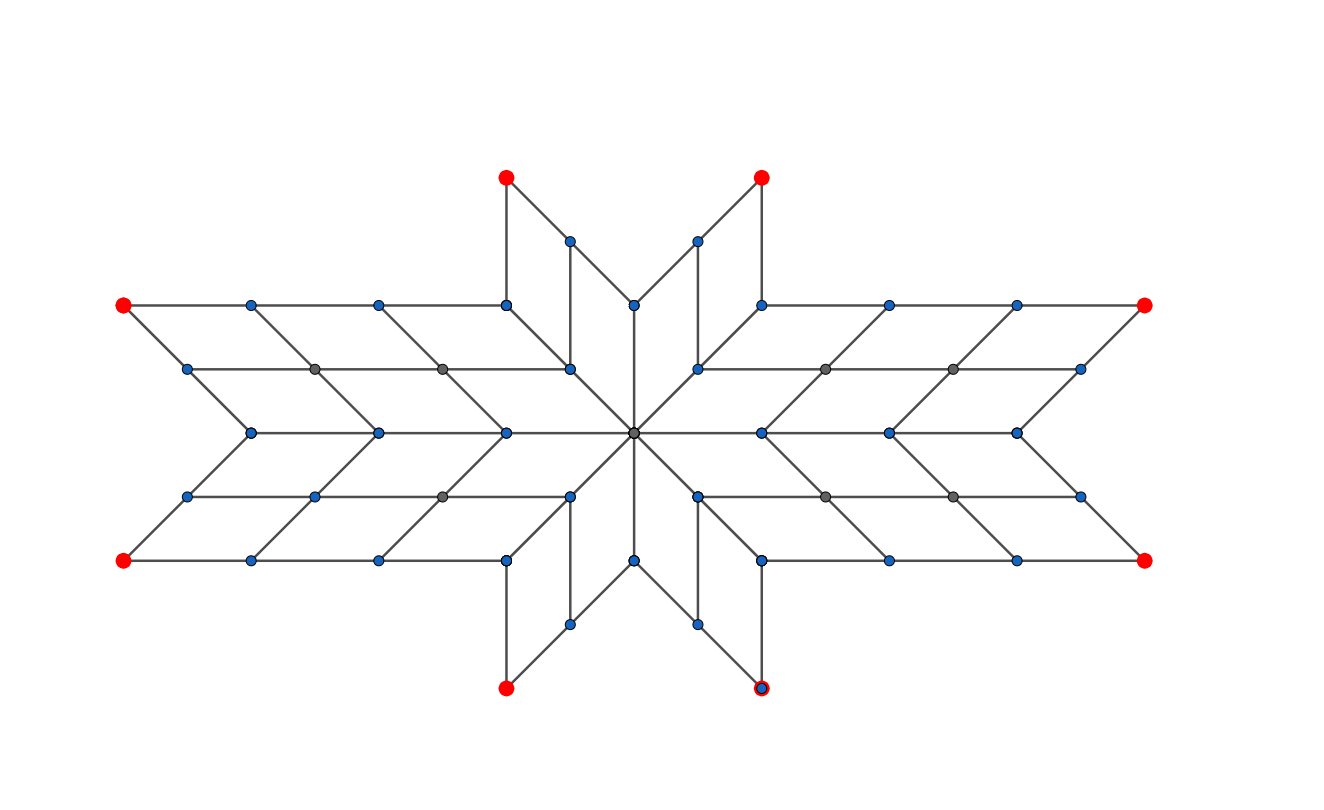}
\end{center}
     \caption{The distance of the singular vertex to the boundary is $1$}
     \label{QUASMESH2}
\end{figure}

\begin{figure}[htbp]
\begin{center}
    \includegraphics[width=0.8\textwidth]{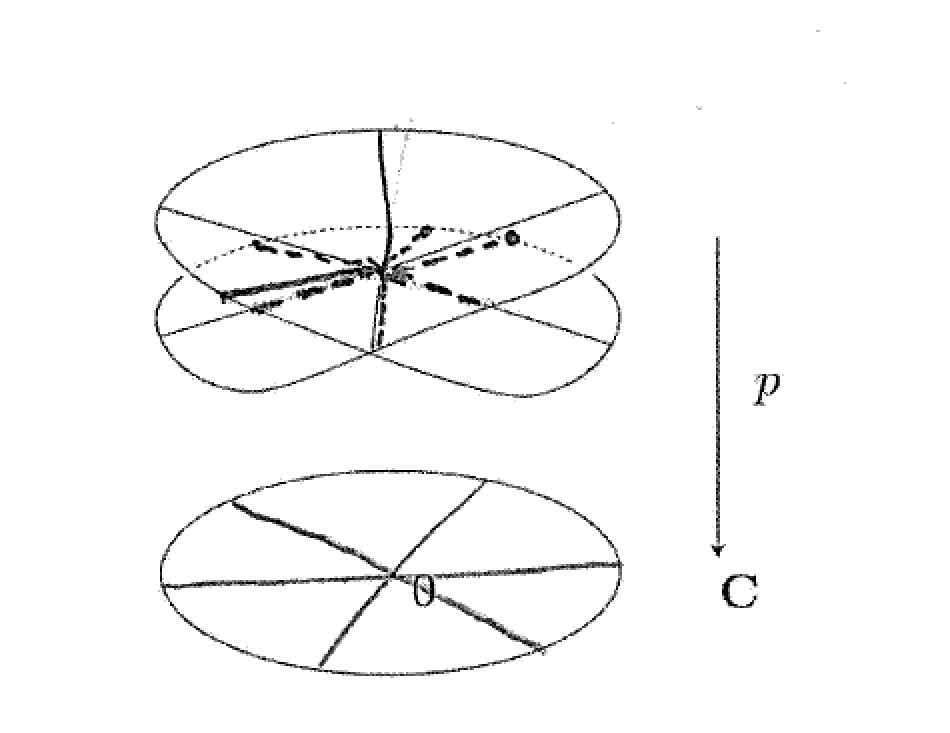}
\end{center}
     \caption{A ramified covering with two sheets. The map $p$ is the projection over the disc whose boundary is the circle $C$}
     \label{RAM}
\end{figure}

\newpage

\bibliographystyle{plain}

\Addresses

\end{document}